\documentclass[11pt]{amsart} 
\usepackage{amsmath}
\usepackage{amssymb}
\usepackage{amsfonts}

\usepackage{verbatim}

\newtheorem{theorem}{Theorem}
\theoremstyle{plain}

\numberwithin{equation}{section}

\input{tcilatex}

\begin{document}
\title[]{Characterization of quadric surfaces in terms of coordinate finite type Gauss map}
\author{Mutaz Al-Sabbagh}
\address{Department of Basic Sciences and Humanities, Imam Abdulrahman bin Faisal University}
\email{malsbbagh@iau.edu.sa}
\author{Hassan Al-Zoubi*}
\address{Department of Mathematics, Al-Zaytoonah University of Jordan, P.O.
Box 130, Amman, Jordan 11733}
\email{dr.hassanz@zuj.edu.jo*}
\date{}

\begin{abstract}
In this article, we introduce an important class of surfaces, namely, quadrics in the Euclidean 3-space $\mathbb{E}^{3}$. We prove that planes, spheres and circular cylinders are the only quadric surfaces whose Gauss map $\boldsymbol{G}$ satisfies a relation of the form $\Delta^{I}\boldsymbol{G}= M \boldsymbol{G}$, where $M$ is a square matrix of order 3 and $\Delta^{I}$ is the Laplace-Beltrami operator corresponding to the first fundamental form $I$ of the surface.
\end{abstract}

\maketitle

\section{Introduction}
Let $\boldsymbol{z}: M^{2}\rightarrow \mathbb{E}^{3},$ be the position vector field of a surface $M^{2}$ in the $3$-dimensional Euclidean space $\mathbb{E}^{3}$. For any two vectors $\boldsymbol{A}=(a_{1},a_{2},a_{3})$ and $\boldsymbol{B}=(b_{1},b_{2},b_{3}) \in \mathbb{E}^{3}$, the inner product on $\mathbb{E}^{3}$ is
\begin{equation}
\boldsymbol{A} \boldsymbol{.} \boldsymbol{B}=  a_{1}b_{1}+a_{2}b_{2}+a_{3}b_{3}. \notag 
\end{equation}
The Euclidean vector product $\boldsymbol{A}\times \boldsymbol{B}$ of $\boldsymbol{A}$ and $\boldsymbol{B}$ is defined as follows:
\begin{equation}
\boldsymbol{A} \times \boldsymbol{B}=  (a_{2}b_{3}-a_{3}b_{2},a_{3}b_{1}-a_{1}b_{3},a_{1}b_{2}-a_{2}b_{1}). \notag
\end{equation}

The concept of surfaces of finite Chen type was born in the year 1973 and became a hot topic of interest in the field of differential geometry and geometric analysis. An Euclidean submanifold is said to be of finite Chen type if its coordinate functions are a finite sum of eigenfunctions of its Laplacian $\Delta$ \cite{C14}. Further, the notion of finite type can be extended to any smooth functions on a submanifold of a Euclidean space or a pseudo-Euclidean space. In this respect, many authors, shed light on the notion of submanifolds of finite type Gauss map. See for example \cite{A21, A7, C9, K2, K3, K4}

Later, a new type of research was generated by investigating surfaces whose Gauss map $\boldsymbol{G}$ satisfies a relation of the form
\begin{equation}  \label{1.1}
\Delta ^{I}\boldsymbol{G}=M \boldsymbol{G},
\end{equation} 
 where $M$ is a square matrix of order 3. In \cite{B4} two classes of surfaces were studied, namely, ruled surfaces and tubes. 

F. Dillen, and others in \cite{D4} studied the class of surfaces of revolution, while in \cite{A18} authors studied the Lorentz-Minkowski version for the same class. Later in \cite{B5, B6, B1} Ch. Baikoussis and L. Verstraelen studied the translation surfaces, the helicoidal surfaces, and the spiral surfaces. In \cite{S0} authors studied translation surfaces of finite type in Sol$_{3}$. H. Al-Zoubi and others investigated the tubes in $\mathbb{E}^{3}$ \cite{A3, A7}. Finally, in \cite{B3} the compact and noncompact cyclides of Dupin were studied.

Following the same ideas of \cite{G3}, it is interesting to study surfaces in $\mathbb{E}^{3}$ whose Gauss map satisfies the relation
\begin{equation}  \label{1.1}
\Delta ^{I}\boldsymbol{G}=M \boldsymbol{G},
\end{equation}%
where $M \in\mathbb{Re}^{3\times3}$. 

In this present paper, we will firstly, create a formula for $\Delta ^{I}\boldsymbol{r}$ and $\Delta ^{I}\boldsymbol{G}$ by using Cartan's
method of the moving frame. Further, we will focus our interest by studyhing the class of quadrics in $\mathbb{E}^{3}$. Our main theorem is 
\begin{theorem}
\label{TH3} Planes, circular cylinders and spheres are the only quadrics in $\mathbb{E}^{3}$ whose Gauss map $\boldsymbol{G}$ satisfying (\ref{1.1}).
\end{theorem}


\section{Basic concepts}

Let $\boldsymbol{r}= \boldsymbol{r}(u,v)$ be a regular parametric representation of a surface $Q$ in $\mathbb{E}^{3}$. A moving frame of the surface $Q$ can be represented by the set of vectors $\hbar =\{\boldsymbol{\xi_{1}}(u,v),\boldsymbol{\xi _{2}}(u,v),\boldsymbol{\xi _{3}}(u,v)\}$, where $\det (\boldsymbol{\xi_{1}},\boldsymbol{\xi _{2}},\boldsymbol{\xi _{3}})=1$. Moreover, we can choose $\boldsymbol{\xi _{3}}$ to be the Gauss map $\boldsymbol{G}$ of $Q$. Hence there exist five linear differential forms $\varpi_{1},\varpi _{2},\varpi _{31},\varpi _{32}$ and $\varpi _{12}$, such that \cite{A5,A2}

\begin{equation*}  \label{2.1}
d\boldsymbol{r}=\varpi _{1}\boldsymbol{\xi _{1}}+\varpi _{2}\boldsymbol{\xi _{2}},\ \ \ d\boldsymbol{G}=\varpi _{31}\boldsymbol{%
\xi _{1}}+\varpi _{32}\boldsymbol{\xi _{2}},
\end{equation*}

\begin{equation*}
d\boldsymbol{\xi _{1}}=\varpi _{12}\boldsymbol{\xi _{2}}-\varpi _{31}\boldsymbol{\xi _{3}},\ \ \ d\boldsymbol{\xi
_{2}}=-\varpi _{12}\boldsymbol{\xi _{1}}-\varpi _{32}\boldsymbol{\xi _{3}},
\end{equation*}

and functions $q_{1},q_{2},a,b,c$ of the variables $u$ and $v$ such that

\begin{equation*}  \label{2.2}
\varpi _{31}=-a\varpi _{1}-b\varpi _{2},\ \ \ \varpi _{32}=-b\varpi_{1}-c\varpi _{2},\ \ \ \varpi _{12}=q_{1}\varpi _{1}+q_{2}\varpi _{2}.
\end{equation*}

We can choose the set $\hbar$, in such way that the principal directions of $Q$ are the vectors $\boldsymbol{\xi_{1}}, \boldsymbol{\xi_{2}}$. Then for the functions $a, c$ and $b$ we get $b=0$ and $a$, $c$ are the principal curvatures of $Q$, hence the differential forms $\varpi_{1}$ and $\varpi_{2}$ reduce to

\begin{equation}  \label{2.3}
\varpi_{1} =-\frac{1}{a}\varpi_{31}, \ \ \ \ \varpi_{2} =-\frac{1}{c}\varpi_{32}.  \notag
\end{equation}

The Gauss curvature and the mean curvature of $Q$ are the following

\begin{equation*}  \label{2.4}
K= ac, \ \ \ \ H=\frac{1}{2}(a+c).
\end{equation*}

We consider a function $h(u,v)\in C^{1}$. Then $\nabla _{1}h,\nabla _{2}h$ denote the derivatives of Pfaff of $h$ along the curves $\varpi _{2}=0,\varpi _{1}=0$ respectively. Thus we have \cite{R2}

\begin{equation}  \label{2.7}
\nabla_{1}\boldsymbol{\xi_{1}}=q_{1}\boldsymbol{\xi_{2}}+a\boldsymbol{G},\ \ \ \nabla_{2}\boldsymbol{\xi_{1}}=q_{2}\boldsymbol{\xi_{2}}+b\boldsymbol{G},
\end{equation}
\begin{equation}  \label{2.8}
\nabla_{1}\boldsymbol{\xi_{2}}=-q_{1}\boldsymbol{\xi_{1}}+b\boldsymbol{G},\ \ \ \ \nabla_{2}\boldsymbol{\xi_{2}}=-q_{2}\boldsymbol{\xi_{1}}+c\boldsymbol{G},
\end{equation}
\begin{equation}  \label{2.6}
\nabla_{1}\boldsymbol{r}=\boldsymbol{\xi_{1}},\ \ \ \nabla_{2}\boldsymbol{r}=\boldsymbol{\xi_{2}},
\end{equation}
\begin{equation}  \label{2.9}
\nabla_{1}\boldsymbol{G}=-a\boldsymbol{\xi_{1}},\ \ \
\nabla_{2}\boldsymbol{G}=-c\boldsymbol{\xi_{2}}.
\end{equation}

The Mainardi-Codazzi equations are

\begin{equation}  \label{2.10}
\nabla_{1}c=q_{2}(a-c),\ \ \ \nabla_{2}a=q_{1}(a-c).
\end{equation}

Let $h$ be a sufficient differentiable function on $Q$. The second Beltrami operator $\Delta^{I}$ of $Q$ is defined by
\begin{equation}  \label{2.11}
\Delta ^{I}h= -\nabla_{1}\nabla_{1}h-\nabla_{2}\nabla_{2}h-q_{2}\nabla_{1}h+q_{1}\nabla_{2}h.
\end{equation}

For the position vector $\boldsymbol{r}$ relation (\ref{2.11}) becomes

\begin{equation}  \label{2.12}
\Delta ^{I}\boldsymbol{r}= -\nabla_{1}\nabla_{1}\boldsymbol{r}-\nabla_{2}\nabla_{2}\boldsymbol{r}-q_{2}\nabla_{1}\boldsymbol{r}+q_{1}\nabla_{2}\boldsymbol{r}.  \notag
\end{equation}

From (\ref{2.6}) we obtain

\begin{equation}  \label{2.13}
\Delta ^{I}\boldsymbol{r}= -\nabla_{1}\boldsymbol{\xi_{1}}-\nabla_{2}\boldsymbol{\xi_{2}}-q_{2}\boldsymbol{\xi_{1}} +q_{1}\boldsymbol{\xi_{2}}.
\end{equation}

Using (\ref{2.7}) and (\ref{2.8}), equation (\ref{2.13}) becomes

\begin{equation*}  \label{2.14}
\Delta ^{I}\boldsymbol{r}=-q_{1}\boldsymbol{\xi_{2}}-a\boldsymbol{G}+q_{2}\boldsymbol{\xi_{1}}- c\boldsymbol{G}-q_{2}\boldsymbol{\xi_{1}}+q_{1}\boldsymbol{\xi_{2}}.
\end{equation*}

Taking into account the last equation and relation (\ref{2.4}), we finally, obtain

\begin{equation}  \label{2.15}
\Delta ^{I}\boldsymbol{r}= -2H\boldsymbol{G}.
\end{equation}


%

We focus our interest now on computing $\Delta ^{I}\boldsymbol{G}$. Inserting the position vector $\boldsymbol{G}$ in (\ref{2.11}) gives

\begin{equation}  \label{2.16}
\Delta ^{I}\boldsymbol{G}= -\nabla_{1}\nabla_{1}\boldsymbol{G}-\nabla_{2}\nabla_{2}\boldsymbol{G}-q_{2}\nabla_{1}\boldsymbol{G}+q_{1}\nabla_{2}\boldsymbol{G}.  \notag
\end{equation}

Using equations (\ref{2.9}), we find
\begin{equation*}  \label{2.17}
\Delta ^{I}\boldsymbol{G}= \nabla_{1}(a\boldsymbol{\xi_{1}})+\nabla_{2}(c\boldsymbol{\xi_{2}})+aq_{2}\boldsymbol{\xi_{1}}-cq_{1}\boldsymbol{\xi_{2}},
\end{equation*}
which becomes
\begin{equation*}  \label{2.18}
\Delta ^{I}\boldsymbol{G}= (\nabla_{1}a)\boldsymbol{\xi_{1}}+a(\nabla_{1}\boldsymbol{\xi_{1}})+(\nabla_{2}c)\boldsymbol{\xi_{2}}+ c(\nabla_{2}\boldsymbol{\xi_{2}})+ aq_{2}\boldsymbol{\xi_{1}}-cq_{1}\boldsymbol{\xi_{2}}.
\end{equation*}

Taking into account equations (\ref{2.7}) and (\ref{2.8}), we get
\begin{equation*}  \label{2.19}
\Delta ^{I}\boldsymbol{G}= (\nabla_{1}a)\boldsymbol{\xi_{1}}+aq_{1}\boldsymbol{\xi_{2}}+a^{2}\boldsymbol{G}+(\nabla_{2}c)\boldsymbol{\xi_{2}}- cq_{2}\boldsymbol{\xi_{1}}+c^{2}\boldsymbol{G}+ aq_{2}\boldsymbol{\xi_{1}}-cq_{1}\boldsymbol{\xi_{2}},
\end{equation*}
or
\begin{equation*}  \label{2.20}
\Delta ^{I}\boldsymbol{G}= (\nabla_{1}a-[c-a]q_{2})\boldsymbol{\xi_{1}}+ (\nabla_{2}c-[c-a]q_{1})\boldsymbol{\xi_{2}} +[a^{2}+c^{2}]\boldsymbol{G}.
\end{equation*}

Using Mainardi-Codazzi equations (\ref{2.10}), we get
\begin{equation*}  \label{2.20}
\Delta ^{I}\boldsymbol{G}= (\nabla_{1}a+\nabla_{1}c)\boldsymbol{\xi_{1}}+ (\nabla_{2}c+\nabla_{2}a)\boldsymbol{\xi_{2}} +(a^{2}+c^{2})\boldsymbol{G}.
\end{equation*}

Once we have
\begin{equation*}  \label{2.21}
(\nabla_{1}a+\nabla_{1}c)\boldsymbol{\xi_{1}}+ (\nabla_{2}c+\nabla_{2}a)\boldsymbol{\xi_{2}} = 2\nabla_{1}H\boldsymbol{\xi_{1}}+ 2\nabla_{2}H\boldsymbol{\xi_{2}} = grad^{I}2H,
\end{equation*}
and
\begin{equation*}  \label{2.22}
4H^{2}-2K = (a^{2}+c^{2}).
\end{equation*}

We finally obtain
\begin{equation}  \label{2.23}
\Delta ^{I}\boldsymbol{G} = grad^{I}2H+(4H^{2}-2K )\boldsymbol{G}.
\end{equation}

\section{Main result}
We consider now a quadric surface $Q$ in $\mathbb{E}^{3}$. Then we have the following three cases 
Case I. $Q$ is ruled, a case that has been studied in  \cite{B2} and it was proved 
\begin{theorem}
\label{TH11} Among the ruled surfaces in $\mathbb{E}^{3}$, the only ones whose Gauss map satisfies (\ref{1.1}) are the planes, and the circular cylinders.
\end{theorem}

Case II. $Q$ is of the form
\begin{equation}  \label{I}
z^{2} = r - p \,x^{2} - q\,y^{2}, \quad p,q,r \in \mathbb{R}, \quad p\,q \neq 0, \quad r > 0,
\end{equation}

Case III. $Q$ is of the form
\begin{equation}  \label{II}
z = \frac{p}{2}\,x^{2} + \frac{q}{2} \,y^{2}, \quad p,q  \in \mathbb{R}, \quad p,q > 0.
\end{equation}

We first prove that a surface of the form (\ref{I}) never satisfies (\ref{1.1}) unless only $p=q=-1$, that is $Q$ is a part of a sphere. Next we prove that a surface of the kind (\ref{II}) is never satisfying (\ref{1.1}). 

\subsection{Quadrics of the first type}

\noindent This type is parameterized as follows

\begin{equation*}
\boldsymbol{r}(u,v)=\left( u,v,\sqrt{r+pu^{2}+qv^{2}}\right).
\end{equation*}

For simplicity, we denote $r+pu^{2}+qv^{2}$ by $\rho $. Then, using the natural frame $\{{\boldsymbol{r}_{u}, \boldsymbol{r}_{v}}\}$ of $Q$ defined by
\begin{equation*}
\boldsymbol{r_{u}}=\left( 1,0,\frac{pu}{\sqrt{\rho}}\right),
\end{equation*}
and
\begin{equation*}
\boldsymbol{r_{v}}=\left( 0,1,\frac{qv}{\sqrt{\rho}} \right),
\end{equation*}
the components $g_{ij}$ of the metric $I$ are 

\begin{equation*}
g_{11}=\boldsymbol{r_{u}} \boldsymbol{.} \boldsymbol{r_{u}}= \frac{\rho +\left( pu\right) ^{2}}{\rho },\ \ \
\end{equation*}
\begin{equation*}
g_{12}=\boldsymbol{r_{u}}\boldsymbol{.} \boldsymbol{r_{v}}=\frac{pquv}{\rho },\ \ \
\end{equation*}
\begin{equation*}
g_{22}=\boldsymbol{r_{v}}\boldsymbol{.} \boldsymbol{r_{v}}=\frac{\rho +\left( qv\right) ^{2}}{\rho }.
\end{equation*}

Hence the Laplacian $\Delta ^{I}$ of $Q$ is \cite{C6}

\begin{eqnarray}  \label{4}
\Delta^{I} &=&-\frac{1}{\Phi}\Bigg[ (\rho + q^{2}v^{2})\frac{\partial ^{2}} {\partial u^{2}} - 2pquv \frac{\partial ^{2}}{\partial u\partial v} +
(\rho + p^{2}u^{2}) \frac{\partial ^{2}} {\partial v^{2}}\Bigg]  \notag \\
&&+\frac{\Omega}{ \Phi^{2}}\Bigg[pu\frac{\partial }{\partial u}+qv\frac{\partial }{\partial v}\Bigg],
\end{eqnarray}
where $\Phi:= \det [g_{ij}]=p(p+1)u^{2}+q(q+1)v^{2}+r$ and $\Omega: = pr + qr + pq(p + 1)u^{2} +pq(q + 1)v^{2}$.

For the normal vector $\boldsymbol{G}$ of $Q$, we have 
\begin{equation*}
\boldsymbol{G}= \frac{\boldsymbol{r_{u}}\times \boldsymbol{r_{v}}}{\sqrt{\Phi}}.
\end{equation*}
After a simple calculations becomes 
\begin{equation*}
\boldsymbol{G}=\frac{1}{\sqrt{\Phi}}\Big(-pu,-qv,\sqrt{\rho}\Big).
\end{equation*}

Let $(n_{1},n_{2},n_{3})$ the components of the vector $\boldsymbol{G}$, and by $\mu _{rs},r,s=1,2,3$ the entries of the matrix $M$. From (\ref{1.1}), we have

\begin{equation}  \label{5}
\Delta ^{I}n_{1}=\Delta ^{I}\Big(-\frac{pu}{\sqrt{\Phi}}\Big)=\mu _{11}\Big(-\frac{pu}{\sqrt{\Phi}}\Big)+\mu _{12}\Big(-\frac{qv}{\sqrt{\Phi}}\Big)+\mu _{13}\Big(\frac{\sqrt{\rho}}{\sqrt{\Phi}}\Big),
\end{equation}

\begin{equation}  \label{6}
\Delta ^{I}n_{2}=\Delta ^{I}\Big(-\frac{qv}{\sqrt{\Phi}}\Big)=\mu _{21}\Big(-\frac{pu}{\sqrt{\Phi}}\Big)+\mu _{22}\Big(-\frac{qv}{\sqrt{\Phi}}\Big)+\mu _{23}\Big(\frac{\sqrt{\rho}}{\sqrt{\Phi}}\Big),
\end{equation}

\begin{equation}  \label{7}
\Delta ^{I}n_{3}=\Delta ^{I}\Big(-\frac{\sqrt{\rho}}{\sqrt{\Phi}}\Big)=\mu _{31}\Big(-\frac{pu}{\sqrt{\Phi}}\Big)+\mu _{32}\Big(-\frac{qv}{\sqrt{\Phi}}\Big)+\mu _{33}\Big(\frac{\sqrt{\rho}}{\sqrt{\Phi}}\Big).
\end{equation}

inserting $n_{1},n_{2}$ of $\boldsymbol{G}$ in (\ref{4}), therefore from (\ref{5}) and (\ref{6}), we conclude

\begin{eqnarray*}
\Delta ^{I}\Big[-\frac{pu}{\sqrt{\Phi}}\Big] &=&-\frac{pu}{\Phi^{\frac{7}{2}}}\left[ p^{2}q(p+1)^{2}(q+1)u^{4}+ f(u,v)\right]  \notag \\
&=&\mu _{11}\Big(-\frac{pu}{\sqrt{\Phi}}\Big)+\mu _{12}\Big(-\frac{qv}{\sqrt{\Phi}}\Big)+\mu_{13}\Big(\frac{\sqrt{\rho}}{\sqrt{\Phi}}\Big),
\end{eqnarray*}
which turns into
\begin{eqnarray}  \label{8}
&&-\frac{pu}{\Phi^{3}}\left[ p^{2}q(p+1)^{2}(q+1)u^{4}+ f(u,v)\right]=  \notag \\
&&-\mu _{11}pu -\mu _{12}qv +\mu_{13}\sqrt{\rho},
\end{eqnarray}
where
\begin{eqnarray}  \label{81}
f(u,v) &=&q^{2}(q+1)^{2}(4p^{2}-3pq+3p-2q)v^{4}  \notag \\
&& +pr(p+1)(2q^{2}+2q+3p+pq)u^{2}  \notag \\
&& +qr(q+1)(-3q^{2}-q+6p+8p^{2}-2pq)v^{2}  \notag \\
&& +pq(p+1)(q+1)(-3q^{2}-q+3p+5pq)u^{2}v^{2}  \notag \\
&& +r^{2}\Big(3p^{2}+3p+q(q+1)+p^{2}+pq\Big).
\end{eqnarray}
\begin{eqnarray*}
\Delta ^{I}\Big(-\frac{qv}{\sqrt{\Phi}}\Big) &=&-\frac{qv}{\Phi^{\frac{7}{2}}}\left[ pq^{2}(p+1)(q+1)^{2}v^{4}+ g(u,v)\right]  \notag \\
&=&\mu _{21}\Big(-\frac{pu}{\sqrt{\Phi}}\Big)+\mu _{22}\Big(-\frac{qv}{\sqrt{\Phi}}\Big)+\mu _{23}\Big(\frac{\sqrt{\rho}}{\sqrt{\Phi}}\Big),
\end{eqnarray*}
which turns into
\begin{eqnarray}  \label{9}
&&-\frac{qv}{\Phi^{3}}\left[ pq^{2}(p+1)(q+1)^{2}v^{4}+ g(u,v)\right]=  \notag \\
&&-\mu _{21}pu-\mu _{22}qv+\mu_{23}\sqrt{\rho},
\end{eqnarray}
where
\begin{eqnarray}  \label{91}
g(u,v) &=&p^{2}(p+1)^{2}(4q^{2}-3pq+3q-2p)u^{4}  \notag \\
&& +qr(q+1)(2p^{2}+2p+3q+pq)v^{2}  \notag \\
&& +pr(p+1)(-3p^{2}-p+6q+8q^{2}-2pq)u^{2}  \notag \\
&& +pq(p+1)(q+1)(-3p^{2}-p+3q+5pq)u^{2}v^{2}  \notag \\
&& +r^{2}\Big(3q^{2}+3q+p(p+1)+qp+q^{2}\Big).
\end{eqnarray}

Using $u=0$ in (\ref{8}), we obtain that
\begin{equation}  \label{10}
-\mu _{12}qv+\mu _{13}\sqrt{1+qv^{2}}= 0.
\end{equation}

Deriving (\ref{10}) with respect to $v$ gives
\begin{equation}  \label{11}
-\mu _{12}\sqrt{1+qv^{2}}+\mu _{13}v = 0.
\end{equation}

Considering (\ref{10}) and (\ref{11}) as a system in $\mu _{12}$ and $\mu _{13}$, and since the determinant
\begin{equation*}
\left|
\begin{array}{cc}
-qv & \sqrt{1+qv^{2}}  \\
-\sqrt{1+qv^{2}} & v
\end{array}%
\right| \neq 0.
\end{equation*}
Therefore we must have $\mu _{12}=\mu _{13}=0$. Hence (\ref{8}) reduces to
\begin{eqnarray*}
\frac{pu}{\Phi^{3}}\left[ p^{2}q(p+1)^{2}(q+1)u^{4}+ f(u,v)\right] =\mu _{11}pu,
\end{eqnarray*}
or
\begin{equation}  \label{12}
 p^{2}q(p+1)^{2}(q+1)u^{4}+ f(u,v) =\mu _{11}{\Phi^{3}}.
\end{equation}

Using $v=0$ in (\ref{12}), and taking into account (\ref{81}), we find
\begin{eqnarray}  \label{13}
\mu _{11}[r+p(p+1)u^{2}]^{3}&=& \Bigg[ p^{2}q(p+1)^{2}(q+1)u^{4} \notag \\
&& +pr(p+1)(2q^{2}+2q+3p+pq)u^{2} \notag \\
&& +r^{2}\Big(3p(p+1)+q(q+1)+p(p+q)\Big)\Bigg].
\end{eqnarray}

Similarly, inserting $v = 0$ in (\ref{9}), gives
\begin{equation*}
-\mu _{21}pu+\mu _{23}\sqrt{1+pu^{2}}= 0.
\end{equation*}

In the same way one can see that $\mu _{21}=\mu _{23}=0$. Then (\ref{9}) turns into
\begin{eqnarray*}
\frac{qv}{\Phi^{3}}\left[ pq^{2}(p+1)(q+1)^{2}v^{4}+ g(u,v)\right]=\mu _{22}qv,
\end{eqnarray*}
or
\begin{equation}  \label{15}
\left[ pq^{2}(p+1)(q+1)^{2}v^{4}+ g(u,v)\right] =\mu _{22}{\Phi^{3}}.
\end{equation}

Inserting $u=0$ in (\ref{15}), and taking into account (\ref{91}), we get
\begin{eqnarray}  \label{16}
\mu _{22}[r+q(q+1)v^{2}]^{3}&=& \Bigg[ pq^{2}(p+1)(q+1)^{2}v^{4} \notag \\
&& +qr(q+1)(2p^{2}+2p+3pq+pq)v^{2} \notag \\
&& +r^{2}\Big(p(p+1)+3q(q+1)+q(p+q)\Big)\Bigg].
\end{eqnarray}

Its clearly that relations (\ref{13}) and (\ref{16}) are polynomials in $u$ and $v$ respectively of degree at most 6. As $p\neq 0, q\neq0$ and $r\neq 0$, then one can be easily obtain that $p= q = -1$. Hence $Q$ is a sphere.

Let $p = q = -1$.  Then from (\ref{81}) and (\ref{91}) respectively, we get $f(u,v)=2r^{2}$ and $g(u,v)=2r^{2}$. Thus from (\ref{12}) and (\ref{15}) we find that $\mu _{11}= \mu _{22}= \frac{2}{r}$.

Besides, relation (\ref{4}) reduces to
\begin{eqnarray}  \label{17}
&&\triangle ^{I}=\frac{1}{r}\Bigg[(u^{2}-r)\frac{\partial ^{2}} {\partial u^{2}} + 2uv \frac{\partial ^{2}}{\partial u\partial v} + \notag \\
&&(v^{2}-r) \frac{\partial ^{2}}{\partial v^{2}}+2u\frac{\partial }{\partial u}+2v\frac{\partial }{\partial v}\Bigg].
\end{eqnarray}
So relation (\ref{7}), becomes

\begin{equation*}  \label{18}
\Delta ^{I}(-\sqrt{\rho})= \frac{2\sqrt{\rho}}{r}=\mu _{31}u+\mu _{32}v+\mu _{33}\sqrt{\rho}.
\end{equation*}

It is easily verified that $\mu _{31}=\mu _{32}=0$ and $\mu _{33}= \frac{2}{r}$. Thus we find that spheres are the only quadric surfaces of the kind (\ref{I}) whose Gauss map satisfies (\ref{1.1}). The resulting matrix is
\begin{equation*}
M =\left[
\begin{array}{ccc}
\frac{2}{r} & 0 & 0 \\
0 & \frac{2}{r} & 0 \\
0 & 0& \frac{2}{r}
\end{array}%
\right].
\end{equation*}

\subsection{Quadrics of the second type}

\noindent This type of surfaces can be parameterized as follows
\begin{equation}  \label{19}
\boldsymbol{r}(u,v) = \left( u,v,\frac{p}{2}u^{2} + \frac{q}{2}v^{2}\right).
\end{equation}

Using the natural frame $\{{\boldsymbol{r}_{u}, \boldsymbol{r}_{v}}\}$ of $Q$ defined by
\begin{equation*}
\boldsymbol{r_{u}}=\left( 1,0,pu\right),
\end{equation*}
and
\begin{equation*} \newline
\boldsymbol{r_{v}}=\left( 0,1,qv\right)
\end{equation*}

\noindent the components $g_{ij}$ of the metric $I$ are 
\begin{equation*}
g_{11}=\left( pu\right) ^{2}+1,\ \ \ g_{12}=pquv,\ \ \ g_{22}=\left(qv\right) ^{2}+1.
\end{equation*}

Therefore the Laplace operator $\Delta^{I}$ of $Q$ is given by
\begin{eqnarray}  \label{20}
\Delta^{I} &=& - \frac{1}{g}\Bigg(Y\frac{\partial ^{2}}{\partial u^{2}} +X\frac{\partial ^{2}}{\partial v^{2}}-2pquv\frac{\partial ^{2}}{\partial u\partial v}\Bigg)  \notag \\
&&+\Bigg[\frac{pY+qX}{g^{2}}\Bigg]\Bigg[pu\frac{\partial }{\partial u}+qv\frac{\partial }{\partial v}\Bigg],
\end{eqnarray}
where
\begin{equation*}
X:=1+p^{2}u^{2},\ \ \ Y:=1+q^{2}v^{2},
\end{equation*}
and
\begin{equation*}
g\colon = \det \left(g_{ij}\right) = p^2 \, u^2 + q^2 \, v^2 +1
\end{equation*}

The Gauss map of $Q$ is
\begin{equation*}
\boldsymbol{G}=\Big({-\frac{pu}{\sqrt{g}},-\frac{qv}{\sqrt{g}},\frac{1}{\sqrt{g}}}\Big).
\end{equation*}

We denote by $(n_{1},n_{2},n_{3})$ the components of $\boldsymbol{G}$, and by $\mu _{rs},r,s=1,2,3$ the entries of the matrix $M$. From (\ref{1.1}), we get

\begin{equation*}  \label{21}
\Delta ^{I}n_{1}=\Delta ^{I}\Big(-\frac{pu}{\sqrt{g}}\Big)=\mu _{11}\Big(-\frac{pu}{\sqrt{g}}\Big)+\mu _{12}\Big(-\frac{qv}{\sqrt{g}}\Big)+\mu _{13}\Big(\frac{1}{\sqrt{g}}\Big),
\end{equation*}

\begin{equation*}  \label{22}
\Delta ^{I}n_{2}=\Delta ^{I}\Big(-\frac{qv}{\sqrt{g}}\Big)=\mu _{21}\Big(-\frac{pu}{\sqrt{g}}\Big)+\mu _{22}\Big(-\frac{qv}{\sqrt{g}}\Big)+\mu _{23}\Big(\frac{1}{\sqrt{g}}\Big),
\end{equation*}

\begin{equation*}  \label{23}
\Delta ^{I}n_{3}=\Delta ^{I}\Big(-\frac{1}{\sqrt{g}}\Big)=\mu _{31}\Big(-\frac{pu}{\sqrt{g}}\Big)+\mu _{32}\Big(-\frac{qv}{\sqrt{g}}\Big)+\mu _{33}\Big(\frac{1}{\sqrt{g}}\Big).
\end{equation*}

Applying the operator $\Delta^{I}$ to the component functions $n_{1}$ and $n_{2}$ of $\boldsymbol{G}$, we find by means of (\ref{20})
\begin{eqnarray*}
\Delta ^{I}\Big(-\frac{pu}{\sqrt{g}}\Big) &=&-\frac{pu}{g^{\frac{7}{2}}}\left[ q^{2}X^{2}+4p^{2}Y^{2}-3q^{3}v^{2}(qX+pY)+5p^{3}q^{3}u^{2}v^{2}+pqg\right]  \notag \\
&=&\mu _{11}\Big(-\frac{pu}{\sqrt{g}}\Big)+\mu _{12}\Big(-\frac{qv}{\sqrt{g}}\Big)+\mu_{13}\Big(\frac{1}{\sqrt{g}}\Big),
\end{eqnarray*}
which turns into
\begin{eqnarray}  \label{24}
\frac{pu}{g^{3}}\left[ q^{2}X^{2}+4p^{2}Y^{2}-3q^{3}v^{2}(qX+pY)+5p^{3}q^{3}u^{2}v^{2}+pqg\right]  \notag \\
=\mu _{11}pu +\mu _{12}qv -\mu_{13},
\end{eqnarray}
and
\begin{eqnarray*}
\Delta ^{I}\Big(-\frac{qv}{\sqrt{g}}\Big) &=&-\frac{qv}{g^{\frac{7}{2}}}\left[ 4q^{2}X^{2}+p^{2}Y^{2}-3p^{3}u^{2}(qX+pY)+5p^{3}q^{3}u^{2}v^{2}+pqg\right]  \notag \\
&=&\mu _{21}\Big(-\frac{pu}{\sqrt{g}}\Big)+\mu _{22}\Big(-\frac{qv}{\sqrt{g}}\Big)+\mu_{23}\Big(\frac{1}{\sqrt{g}}\Big),
\end{eqnarray*}
which turns into
\begin{eqnarray}  \label{25}
\frac{qv}{g^{3}}\left[ 4q^{2}X^{2}+p^{2}Y^{2}-3p^{3}u^{2}(qX+pY)+5p^{3}q^{3}u^{2}v^{2}+pqg\right]  \notag \\
=\mu _{21}pu +\mu _{22}qv -\mu_{23}.
\end{eqnarray}

Inserting $u=0$ in (\ref{24}), then the left side of the equation (\ref{24}) vanishes. Therefore we are left to
\begin{equation*}
\mu _{12}qv-\mu _{13}= 0,
\end{equation*}
which immplies that $\mu _{12}=\mu _{13}= 0$. So equation (\ref{24}) becomes
\begin{eqnarray}  \label{26}
&& q^{2}X^{2}+4p^{2}Y^{2}-3q^{3}v^{2}(qX+pY)  \notag \\
&& +5p^{3}q^{3}u^{2}v^{2}+pqg-\mu _{11}g^{3}=0.
\end{eqnarray}

Similarly, if we put $v=0$ in (\ref{25}), then the left side of (\ref{25}) vanishes. In the same way equation (\ref{25}) becomes
\begin{eqnarray}  \label{27}
&& 4q^{2}X^{2}+p^{2}Y^{2}-3p^{3}u^{2}(qX+pY)  \notag \\
&& +5p^{3}q^{3}u^{2}v^{2}+pqg-\mu _{22}g^{3}=0.
\end{eqnarray}

Equations (\ref{26}) and (\ref{27}) are nontrivial polynomials in $u$ and $v$ with constant coefficients. These two polynomials can never be zero, unless $p =q= 0$, which is clearly impossible since $p,q > 0$.

\section{Conclusion}
This research article was divided into three sections, where after the introduction, the needed definitions and relations regarding this interesting field of study were given. Then a formula for the Laplace operator corresponding to the first fundamental form $I$ was proved once for the position vector and another for the Gauss map of a surface $Q$ by using Cartan's method of the moving frame. Finally, we classify the quadric surfaces $Q$ satisfying the relation $\Delta \boldsymbol{G}= M\boldsymbol{G}$, for a real square matrix $M$ of order 3. An interesting study can be drawn, if this type of study can be applied to other classes of surfaces that have not been investigated yet such as spiral surfaces, or tubular surfaces.


\end{document}